\documentclass[12pt]{article}
\usepackage[dvips]{color}
\usepackage{graphicx}
\usepackage{epsfig}
\usepackage{verbatim}
\usepackage{amsthm}
\usepackage{amssymb}
\usepackage{amsxtra}     
\usepackage{multirow}

\newtheorem{theorem}{Theorem}[section]
\newtheorem{definition}[theorem]{Definition}

\newtheorem{corollary}[theorem]{Corollary}

\newcommand{\PROOF}{\noindent {\bf Proof}: }

\DeclareSymbolFont{AMSb}{U}{msb}{m}{n}
\DeclareMathSymbol{\N}{\mathbin}{AMSb}{"4E}
\DeclareMathSymbol{\Z}{\mathbin}{AMSb}{"5A}
\DeclareMathSymbol{\R}{\mathbin}{AMSb}{"52}
\DeclareMathSymbol{\Q}{\mathbin}{AMSb}{"51}
\DeclareMathSymbol{\I}{\mathbin}{AMSb}{"49}
\DeclareMathSymbol{\C}{\mathbin}{AMSb}{"43}

\begin{document}

\title{A Note on the Hodge Structure of the Intersection of Coloring Complexes}
\author{Sarah Crown Rundell \\ Denison University} \maketitle

\begin{abstract}
\hspace{.5in}Let $G$ be a simple graph with $n$ vertices.  The coloring complex $\Delta(G)$ was defined by Steingr\'{\i}msson, and the homology of $\Delta(G)$ was shown to be nonzero only in dimension $n-3$ by Jonsson.  Hanlon recently showed that the Eulerian idempotents provide a decomposition of the homology group $H_{n-3}(\Delta(G))$ where the dimension of the $j^{th}$ component in the decomposition, $H_{n-3}^{(j)}(\Delta(G))$, equals the absolute value of the coefficient of $\lambda^{j}$ in the chromatic polynomial of $G$, $\chi_{G}(\lambda)$.

\hspace{.5in}Jonsson recently studied the topology of intersections of coloring complexes.  In this note, we show that the coefficient of the ${j}^{th}$ term in the chromatic polynomial of the intersection of coloring complexes gives the Euler Characteristic of the $j^{th}$ Hodge subcomplex of the Hodge decomposition of the intersection of coloring complexes. 
\end{abstract}

\section{Introduction}

\hspace{.5in}Consider $R=A/I$ where $A=F[x_{S}|S \subseteq \{
1,...,n \} ]$, $I$ is the ideal generated by $\{x_{U}x_{T}|U \not
\subseteq T, T \not \subseteq U \} $, and $F$ is a field of
characteristic zero. Then consider the ideal $K_G$ generated by the
monomials $x_{X_{1}}^{e_{1}}x_{X_{2}}^{e_{2}}...x_{X_{l}}^{e_{l}}$,
$e_{i} > 0$, such that for all $i$, $1 \leq i \leq l+1$, $Y_{i} =
X_{i}\backslash X_{i-1}$ does not contain an edge of $G$ ($X_{0} = \emptyset$ and $X_{l+1} = \{1,...,n\}$).
Steingr\'{\i}msson~\cite{sm} calls $K_{G}$ the coloring ideal, since
there is a bijection between monomials of $K_{G}$ of degree $r$ and
colorings of $G$ with $r+1$ colors.  He then notes that the quotient
$R/K_{G}$ is the face ring of a simplicial complex, $\Delta(G)$.  This complex is called the coloring complex of $G$.

\hspace{.5in}In 2005, Jonsson~\cite{jo}, showed that $\Delta(G)$ has the homology of a wedge of $a_{G}-1$ spheres, where $a_{G}$ is the number of acyclic orientations of $G$.
Stanley~\cite{st} showed that $a_G$ is $(-1)^n\chi_{G}(-1)$, and Jonsson~\cite{jo} concluded from
this result that the dimension of the $(n-3)^{rd}$ homology group of
$\Delta(G)$ is in fact, $(-1)^n\chi_{G}(-1)-1$.  Hanlon~\cite{ha} showed that Eulerian idempotents provide a decomposition of the homology group $H_{n-3}(\Delta(G))$ where the dimension of the $j^{th}$ component in the decomposition, $H_{n-3}^{(j)}(\Delta(G))$, equals the absolute value of the coefficient of $\lambda^{j}$ in the chromatic polynomial of $G$, $\chi_{G}(\lambda)$.

\hspace{.5in}Recently, Jonsson~\cite{j2} studied the topology of the intersections of coloring complexes and determined a condition for which these complexes are homology Cohen-Macaulay.  In this note, we provide a generalization of Hanlon's result to the intersection of coloring complexes.

\section{The Coloring Complex}

\hspace{.5in}We begin by
defining Steingr\'{\i}msson's~\cite{sm} coloring complex following
the presentation in Jonsson~\cite{jo}.

\hspace{.5in}Let $G$ be a simple graph on $n$ vertices.  Let $ (B_{1},...,B_{r+2}) $ be an ordered partition
of $ \{ 1,...,n \} $ where at least one of the $B_{i}$ contains an edge of $G$, and let $\Delta_{r}(G)$ be the set of ordered partitions
$ ( B_{1},...,B_{r+2} ) $.

\begin{definition} The \emph{coloring complex} of $G$, denoted $\Delta(G)$, is the simplicial complex defined by the sequence:

\begin{center}
$...\rightarrow C_{r} \stackrel{\delta_{r}}{\rightarrow} C_{r-1}
\stackrel{\delta_{r-1}}{\rightarrow} ...
\stackrel{\delta_1}{\rightarrow} C_{0}
\stackrel{\delta_0}{\rightarrow} C_{-1}
\stackrel{\delta_{-1}}{\rightarrow} 0$
\end{center}

where $C_{r}$ is the vector space over a field of characteristic
zero with basis $\Delta_{r}(G)$ and

\begin{center}
$\partial_r( (B_{1},...,B_{r+2}) ) := \displaystyle
\sum_{i=1}^{r+1} (-1)^i ( B_1,...,B_{i} \bigcup B_{i+1},...,
B_{r+2} ) $.
\end{center}

\end{definition}

Notice that $\partial_{r-1} \circ \partial_{r} = 0$.  Then:

\begin{definition}
The $r^{th}$ homology group of $\Delta(G)$, $H_{r}(\Delta(G)) :=
ker(\partial_{r})/im(\partial_{r+1})$.
\end{definition}

The intersection of coloring complexes is defined as follows:

\begin{definition}  Let $m \in \N$ and let $\mathbf{G} = (G_1, ..., G_m)$ be a sequence of nonempty graphs on a vertex set, $V$, of cardinality $n$.  The \emph{intersection of coloring complexes} $ of \mathbf{G}$, denoted $\Delta(\mathbf{G})$, is 
$$\Delta(\mathbf{G}) = \bigcap_{i=1}^{m} \Delta(G_i).$$
\end{definition}

\hspace{.5in}Jonsson~\cite{j2} showed that when the sequence $\mathbf{G}$ is diagonally cycle-free, the complex $\Delta(\mathbf{G})$ is homology Cohen-Macaulay.

\begin{definition}  Let $\mathbf{G} = (G_1, ..., G_m)$ be a sequence of nonempty graphs on a vertex set $V$, and let $E = \{ e_1, ..., e_m \}$ be a set of edges where $e_i \in G_i$ for each $i$, $1 \leq i \leq m$.  The set $E$ is called a \emph{diagonal} of $\mathbf{G}$.  The sequence $\mathbf{G}$ is \emph{diagonally cycle-free} if the edge sets of the $G_i$ are mutually disjoint, and if for all $E$, the graph $(V, E)$ is acyclic.
\end{definition}

\hspace{.5in}The chromatic polynomial, $\chi_{\mathbf{G}}(\lambda)$, of $\mathbf{G}$ is defined to be the number of ways of coloring the vertices in $V$ with $\lambda$ colors such that the coloring is proper for at least one of the $G_i$.  Let $\mathbf{G'} = (G_1, ..., G_m)$.  It is straightforward to prove the following recursive identity:
$$\chi_{\mathbf{G}}(\lambda) = \chi_{\mathbf{G'}, G_{m-1}}(\lambda) - \chi_{\mathbf{G'}, G_{m-1} \cup G_m}(\lambda) + \chi_{\mathbf{G'}, G_m}(\lambda).$$

\section{Eulerian Idempotents}

\hspace{.5in} Recently, Hanlon~\cite{ha} showed that there is a Hodge decomposition of $H_{n-3}(\Delta(G))$ for a graph $G$; we will discuss this result and its generalization to the intersection of coloring complexes.  In order to describe this decomposition, we must first define and describe the Eulerian idempotents.  The Eulerian idempotents have many interesting properties and have proved useful in many different algebraic and combinatorial problems.  For more information on Eulerian idempotents see ~\cite{gs}, ~\cite{lo}, ~\cite{h1}, and ~\cite{h2}.

\hspace{.5in}Define a \emph{descent} of a permutation $\pi \in S_{n}$ to be a couple of consecutive numbers $(i,i+1)$ such that $\pi(i)>\pi(i+1)$. It follows from Loday's ~\cite{lo} definition that the Eulerian idempotents $e_r^{(j)}$ can be defined by the identity:

\begin{definition} The Eulerian idempotents are defined by 
$$\sum_{j=1}^n t^j e_n^{(j)} =
\sum_{\pi\in S_n} \left(\stackrel{n+t-des(\pi)-1}{n}\right)\text{sgn}(\pi)\pi,$$ where \text{des}$(\pi)$ is the number of descents of $\pi$.
\end{definition}

\hspace{.5in}There are several important properties of the Eulerian idempotents which are due to Gerstenhaber and Schack~\cite{gs}.  In their paper, they show that the Eulerian idempotents are mutually orthogonal idempotents and that their sum is the unit element in $\C[S_n]$.  So then for any $S_{n}$-module, $M$, we have that

\begin{center} $\displaystyle M = \bigoplus_{j} e_{n}^{(j)}M$. \end{center}

\hspace{.5in}Notice that we can define an action of $S_{r+2}$ on $\Delta_{r}(G)$.  Namely, if $\sigma \in S_{r+2}$, then $\sigma \cdot ( B_{1},...,B_{r+2} ) = ( B_{\sigma^{-1}(1)},...,B_{\sigma^{-1}(r+2)} ) $.  This action then makes $C_{r}$ into an $S_{r+2}$-module.

Hanlon~\cite{ha} notes (and this result can be derived from the work of Gerstenhaber and Schack~\cite{gs}) in Lemma $2.1$ of his paper that for any graph $G$ and for each $r,j$,

\begin{center} $\partial_{r} \circ e_{r+2}^{(j)} = e_{r+1}^{(j)} \circ \partial_{r}$. \end{center}

This implies then that, for each $j$, $C_{r}^{(j)}=e_{r+2}^{(j)}C_{r}$ is a subcomplex of $(C_{*}(\Delta(G)),\partial_{*})$.  We may then consider the homology of the subcomplex, and it will be denoted by $H_{*}^{(j)}(\Delta(G))$.  So then we have
$$H_{r}(\Delta(G)) = \oplus_{j}  H_{r}^{(j)}(\Delta(G)).$$

The above decomposition is called the \emph{Hodge decomposition} of $H_{*}(\Delta(G))$.

\hspace{.5in}Hanlon~\cite{ha} showed that there is a Hodge decomposition of the top homology group of $\Delta(G)$, i.e. $\displaystyle H_{n-3}(\Delta(G)) = \oplus_{j=1}^{n-1} H_{n-3}^{(j)}(\Delta(G))$.  Further, he showed that the dimension of the $j^{th}$ Hodge piece is equal to the absolute value of the coefficient of $\lambda^{j}$ in the chromatic polynomial of $G$.

\hspace{.5in}In the case of the intersection of coloring complexes, the Hodge decomposition of $\Delta(\mathbf{G})$ can be defined by following the same process as above.  For an arbitrary sequence of nonempty graphs, $\mathbf{G} = (G_1, ..., G_m)$, the homology of $\Delta(\mathbf{G})$ is not concentrated in one dimension, and thus Hanlon's result does not in general hold.  In this note, we will provide a generalization of Hanlon's result.  We also will note that when the sequence $\mathbf{G}$ is diagonally cycle-free the absolute value of the coefficient of $\lambda^j$ in $\chi_{\mathbf{G}}(\lambda)$ is equal to the dimension of $H_{*}^{(j)}(\Delta(\mathbf{G}))$.

\section{The Relationship Between the Chromatic Polynomial of $\Delta(\mathbf{G})$ and $\Delta(\mathbf{G})$}

\hspace{.5in}In this section, we will provide a generalization of Hanlon's result to the intersection of coloring complexes.  In our study, we will need the following definition:

Let $X^{(j)}(\Delta(\mathbf{G}))$ denote the Euler Characteristic of the $j^{th}$ Hodge piece of $\Delta(\mathbf{G})$.  In particular,
\begin{eqnarray*}
X^{(j)} & =  & \sum_{i=-1}^{n-r-1} (-1)^{i} \dim(C_{i}^{(j)}(\Delta(\mathbf{G})) \\
& = & \sum_{i=-1}^{n-r-1} (-1)^{i} \dim(H_{i}^{(j)}(\Delta(\mathbf{G}))
\end{eqnarray*}

\begin{theorem} \label{EulerChar}
For each $j$,
\begin{center}$X^{(j)}(\Delta(\bf{G})) = - [ \lambda^{j} ]( \chi_{\mathbf{G}}(-\lambda)-(-\lambda)^n) $. \end{center}
where $ [ \lambda^{j} ] (\chi_{\mathbf{G}}(-\lambda) - (-\lambda)^n) $ denotes the coefficient of $\lambda^{j}$ in $\chi_{\mathbf{G}}(-\lambda) - (-\lambda)^n$.
\end{theorem}

\PROOF
The proof is similar to the proof of Theorem 4.1 in ~\cite{ha} and will be by induction on $m$.  The base case follows from Theorem 4.1 in ~\cite{ha}.  By way of induction, suppose $m \geq 2$.  Let $\mathbf{G'} = (G_1, ..., G_{m-2})$.  Jonsson~\cite{j2} notes: 
$$\Delta_r(\mathbf{G'}, G_{m-1} \cup G_m) =  (\Delta_r(\mathbf{G'}, G_{m-1}) \backslash \Delta_r(\mathbf{G})) \cup \Delta_r(\mathbf{G'},G_m) .$$

Let $C_{r}(\mathbf{G})$ be the vectorspace having as its basis $\Delta_r(\mathbf{G})$.  This then implies that
$$C_r(\mathbf{G'}, G_{m-1} \cup G_m) =  (C_r(\mathbf{G'}, G_{m-1}) \backslash C_r(\mathbf{G})) \oplus C_r(\mathbf{G'},G_m) .$$
Notice that this isomorphism commutes with the action of $S_{r+2}$, which implies that
$$C_{r}^{(j)}(\mathbf{G'}, G_{m-1} \cup G_m) =  (C_{r}^{(j)}(\mathbf{G'}, G_{m-1}) \backslash C_{r}^{(j)}(\mathbf{G})) \oplus C_{r}^{(j)}(\mathbf{G'},G_m) .$$

Therefore,
$$\dim C_{r}^{(j)}(\mathbf{G'}, G_{m-1} \cup G_m) =  \dim C_{r}^{(j)}(\mathbf{G'}, G_{m-1}) - \dim C_{r}^{(j)}(\mathbf{G}) + \dim C_{r}^{(j)}(\mathbf{G'},G_m) .$$

In particular, 
$$\dim C_{r}^{(j)}(\mathbf{G}) = \dim C_{r}^{(j)}(\mathbf{G'}, G_{m-1}) - \dim C_{r}^{(j)}(\mathbf{G'}, G_{m-1} \cup G_m) + \dim C_{r}^{(j)}(\mathbf{G'},G_m) .$$

Then,
\begin{eqnarray*}
X^{(j)}(\Delta(\mathbf{G})) & = & \sum_{i=-1}^{n-3} (-1)^{i} \dim(H_{i}^{(j)}(\Delta(\mathbf{G})) \\
& = &  \sum_{i=-1}^{n-3} (-1)^{i}  \dim(C_{i}^{(j)}(\Delta(\mathbf{G})) \\
& =  & \sum_{i=-1}^{n-3} (-1)^{i}(  \dim C_{r}^{(j)}(\mathbf{G'}, G_{m-1}) - \dim C_{r}^{(j)}(\mathbf{G'}, G_{m-1} \cup G_m) + \dim C_{r}^{(j)}(\mathbf{G'},G_m)) \\
& =  & \sum_{i=-1}^{n-3} (-1)^{i}(  \dim C_{r}^{(j)}(\mathbf{G'}, G_{m-1}) - \sum_{i=-1}^{n-3} (-1)^{i}( \dim C_{r}^{(j)}(\mathbf{G'}, G_{m-1} \cup G_m)) \\
& & + \sum_{i=-1}^{n-3} (-1)^{i}(\dim C_{r}^{(j)}(\mathbf{G'},G_m)) \\
& = & -[\lambda^{j}] ((\chi_{\mathbf{G'}, G_{m-1}}(-\lambda) - (-\lambda)^n) - (\chi_{\mathbf{G'}, G_{m-1} \cup G_m}(-\lambda)-(-\lambda)^n) + (\chi_{\mathbf{G'}, G_m}(-\lambda) - (-\lambda)^n))\\
& = & -[\lambda^{j}] (\chi_{\mathbf{G}}(-\lambda) -(-\lambda)^n).
\end{eqnarray*}

where the second to last inequality follows by induction and the last line follows by identity (12) in Jonsson~\cite{j2}.

\qed

\begin{corollary}  If $\mathbf{G}$ is diagonally cycle free, then for each $j$, the dimension of $H_{n-3}^{(j)}(\Delta(\mathbf{G}))$ equals the absolute value of $\lambda^{j}$ in $\chi_{\mathbf{G}}(-\lambda)$.
\end{corollary}

\bibliography{biblio}

\end{document}